\documentclass{article}

\usepackage{amsmath}
\usepackage{amsthm}
\usepackage{amssymb}
\usepackage{hyperref} 
\usepackage{srcltx}
\usepackage{color}

\def\Sch{{\hbox{\tiny\rm Sch}}}
\def\cI{{\cal I}}
\def\cS{{\cal S}}

\def\inter{{\hbox{\rm int}\,}}

\def\cO{{\cal O}}
\def\cM{{\cal M}}

\def\bR{{\mathbf{R}}}
\def\Dom{{\hbox{\rm Dom}\,}}

\def\bE{{\mathbf{E}}}

\def\cN{{\cal N}}

\def\cF{{\cal F}}
\def\cM{{\cal M}}

\def\bE{{\mathbf{E}}}
\def\bR{{\mathbf{R}}}
\def\cC{{\cal C}}
\def\cP{{\cal P}}
\def\Risk{{\hbox{\rm Risk}}}
\def\thRisk{{\mbox{\em Risk}}}
\def\Opt{{\hbox{\rm Opt}}}
\def\thOpt{{\mbox{\em Opt}}}

\def\sign{{\hbox{\rm sign}}}

\newtheorem{propos}{Proposition}
\newtheorem{cor}{Corollary}
\newtheorem{rem}{Remark}

\title{On Lower Complexity Bounds for Large-Scale Smooth Convex Optimization\thanks{The research was supported by the NSF grant  CMMI-1232623.}}
\author{Crist\'obal Guzm\'an\thanks{cguzman@gatech.edu} ~~and~
	    Arkadi Nemirovski\thanks{arkadi.nemirovski@isye.gatech.edu}\\
	    \small{H. Milton Stewart School of Industrial and Systems Engineering,}\\
	    \small{Georgia Institute of Technology, Atlanta, GA, USA.}\\}

\date{}

\begin{document}
\maketitle

\begin{abstract}
We derive lower bounds on the black-box oracle
complexity of large-scale smooth convex minimization
problems, with emphasis on minimizing smooth (with H\"older continuous, with a given exponent and constant, gradient) convex functions over high-dimensional $\|\cdot\|_p$-balls, $1\leq p\leq \infty$. Our bounds turn out to be tight (up to logarithmic in the design dimension factors), and can be viewed as a substantial extension of the existing
lower complexity bounds for large-scale convex minimization covering the nonsmooth case and the ``Euclidean'' smooth case (minimization of convex functions with Lipschitz continuous gradients over Euclidean balls). As a byproduct of our results, we demonstrate that the classical Conditional Gradient algorithm is near-optimal, in the sense of Information-Based Complexity Theory, when minimizing smooth convex functions over high-dimensional $\|\cdot\|_\infty$-balls and their matrix analogies -- spectral norm balls in the spaces of square matrices.
\end{abstract}

\section{Introduction}
Huge sizes of convex optimization problems arising in some modern applications (primarily, in big-data-oriented signal processing and machine learning) are beyond
the ``practical grasp'' of the state-of-the-art Interior Point Polynomial Time methods with their computationally demanding iterations. Indeed, aside of rare cases of problems with ``extremely favourable'' structure, the arithmetic cost of an interior point iteration is at least cubic in the design dimension $n$ of the instance; with $n$ in the range of $10^4$ -- $10^6$, as is the case in the outlined applications, this makes a single iteration ``lasting forever.'' The standard techniques for handling large-scale convex problems -- those beyond the practical grasp of Interior Point methods -- are First Order methods (FOM's). Under favorable circumstances, iterations of FOM's are much cheaper than those of interior point methods, and the convergence rate, although just sublinear, is fully or nearly dimension-independent, which makes FOM's the methods of choice when medium-accuracy solutions
to large-scale convex programs are sought. Now, as a matter of fact, all known FOM's are ``black-box-oriented'' -- they ``learn'' the problem being solved solely via the local information (values and (sub)gradients of the objective and the constraints) accumulated along the search points generated by the algorithm. As a result, ``limits of performance'' of FOM's are governed by Information-Based Complexity Theory.
Some basic results in this direction have been
established in the literature \cite{Nemirovski:1983}; in particular, we know well enough what is the Information-Based Complexity of natural families of
convex minimization problems $\min_{x\in X} f(x)$ with {\sl nonsmooth} Lipschitz continuous objectives $f$ and how the complexity depends on the geometry and the dimension of
the domain $X$. In the smooth case, our understanding is somehow limited; essentially, tight {\sl lower} complexity bounds are known only in the case when $X$ is Euclidean ball and $f$ is convex function with Lipschitz continuous gradient. Lower bounds here come from least-squares problems \cite{Nemirovski:1991,Nemirovski:1992}, and the underlying techniques for generating ``hard instances'' heavily utilize the rotational invariance of a Euclidean ball.
\par
In this paper, we derive tight lower bounds on information-based complexity of families of convex minimization problems $\{\min_{x\in X} f(x): f\in\cF\}$, where $X$ is $n$-dimensional $\|\cdot\|_p$-ball, $1\leq p\leq\infty$, and $\cF$ is the family of all continuously differentiable convex objectives with given smoothness parameters (H\"older exponent and constant). We believe that these bounds could be of interest in some modern applications, like $\ell_1$  and nuclear norm minimization in Compressed Sensing, where one seeks to minimize a smooth, most notably, quadratic convex function over high-dimensional $\ell_1$-ball in $\bR^n$ or nuclear norm ball in the space of $n\times n$ matrices. Another instructive application of our results is establishing the near-optimality, in the sense of information-based complexity, of Conditional Gradient (a.k.a. Frank-Wolfe) algorithm as applied to minimizing smooth convex functions over large-scale boxes (or unit balls of spectral norm on the space of matrices)\footnote{Originating from \cite{Frank:1956}, Conditional Gradient algorithm was intensively studied in 1970's (see \cite{Dem:Rub:1970,Pshe:1994} and references therein); recently, there is a significant burst of interest in this technique, due to its ability to handle smooth large-scale convex programs on ``difficult geometry'' domains, see \cite{Hazan:2008,Jaggi:2011,Jaggi:2013,HJN:2013,CJN:2013} and references therein.}.

\subsection{Contributions}

Our first contribution is a unified framework to prove lower
bounds for a variety of domains and different smoothness
parameters of the objective with respect to a norm (for
consistency we use the norm induced by the domain).
In order to construct hard instances for lower bounds
we need the normed space under consideration to satisfy
a ``smoothing property.'' Namely, we need the existence of a
``smoothing kernel'' -- a convex function with Lipschitz continuous gradient and ``fast growth.''
These properties guarantee that the inf-convolution
\cite{HiriartUrruty:2001}
of a Lipschitz continuous convex function $f$ and the smoothing
kernel is smooth, and its local behaviour depends only on
the local behavior of $f$. A novelty here, if any, stems from the
fact that we need Lipschitz continuity of the gradient w.r.t.
a given, not necessarily Euclidean, norm, while the standard
Moreau envelope technique is adjusted to the case of the
Euclidean norm\footnote{$^)$It well may happen that the
extensions of the classical Moreau results  which we present
in Section \ref{sec:Smoothing}   are known, so that the material
in this section does not pretend to be novel. This being said,
at this point in time we do not have at our disposal references
to the results on smoothing we need, and  therefore we decided
to augment these simple results with their proofs, in order to
make our presentation self-contained.}$^)$

We establish lower bounds on complexity of smooth convex minimization for general spaces satisfying the
smoothing property. Our proof mimics the construction of hard instances for
nonsmooth convex minimization \cite{Nemirovski:1983}, which now  are
properly smoothed by the inf-convolution.

With this general result, we are able to provide a unified
analysis for lower bounds for smooth convex minimization
over $n$-dimensional $\|\cdot\|_p$-balls, $1\leq p \leq\infty$.
We show that in the large-scale case, our lower complexity bounds match, within at worst a  logarithmic in $n$ factor,
the upper complexity bounds associated with Nesterov's fast gradient algorithms \cite{Nemirovski:1985,Elster:1993}. When
$p=\infty$, this result implies near optimality
of the Conditional Gradient  algorithm.

As a final application, we point out how our lower bounds extend
to matrix optimization under Schatten norm constraints.

\subsection{Related  work}

\noindent{\bf Oracle Complexity:}
The analysis of convex optimization algorithms
via oracle complexity and lower complexity
bounds were first studied in \cite{Nemirovski:1983}.
Other standard references are \cite{Nemirovski:1994,
Nesterov:2004}. The oracle complexity of smooth
convex optimization over Euclidean domains was
studied in \cite{Nemirovski:1983,
Nemirovski:1991,Nemirovski:1992}.

For optimal methods under non-Euclidean domains for
smooth spaces and $p$-norms, where $2\leq p<\infty$
we refer to \cite{Elster:1993}
(for the case $p=2$ there is an interesting new algorithm
that adapts itself to the smoothness parameter in the
objective \cite{Nesterov:2013}).

It should be mentioned that for the case $2\leq p<\infty$
the lower bounds in this paper were announced
in \cite{Nemirovski:1985,Elster:1993} (and proved
by the second author of this paper);
however, aside of the very special case of $p=2$, the
highly technical original  proofs of the bounds were
never published. For this reason,
we recently have revisited the original proofs and
were able to simplify them dramatically, thus
making them publishable.\\

\noindent{\bf The Conditional Gradient algorithm and
complexity under Linear Optimization oracles:}
The recent body of work on the Conditional Gradient
algorithm is enormous. For upper bounds on its complexity
we refer to \cite{Clarkson:2008,Hazan:2008,Jaggi:2013,
Lan:2013, Garber:2013}. Interestingly, the last two
references include results on linear convergence of
the Conditional Gradient method for the strongly convex case,
accelerated methods based on Linear Optimization
oracles, and applications
to stochastic and online convex programming.

Besides these accuracy upper bounds, there are some
interesting lower bounds for algorithms based on a Linear
Optimization oracle (whose only assumption is that the
Linear Optimization oracle returns a solution that is a vertex of the
domain): some of these contributions can be found in
\cite{Jaggi:2013,Lan:2013}. Observe that a Linear Optimization oracle
is in general less powerful than an arbitrary local oracle
(in particular the first-order one) considered in our paper, and thus their lower
bounds do not imply ours. However, our result for $p=\infty$
improves on their lower bounds (disregarding logarithmic
factors).
\\

\subsection{Notation and preliminaries} \label{prelims}

\noindent{\bf Algorithms and Complexity:}
In the black-box oracle complexity model for convex optimization
we are interested in solving problems of the form
$$
\Opt(f)=\min_{x\in X} f(x)\eqno{(P_{f,X})}
$$
where $X$ is a given convex compact subset of a normed space $(\bE,\|\cdot\|)$,
and $f$ is known to belong to a given family $\cF$ of continuous convex
functions on $\bE$. This defines the family of problems ${\cal P}({\cal F},X)$
comprised of problems $(P_{f,X})$ with $f\in \cF$.
We assume that the  family $\cF$  is equipped with an
{\sl oracle} $\cO$ which, formally, is a function $\cO(f,x)$ of $f\in \cF$
and $x\in \bE$ taking values in some {\sl information space} $\cI$;
when solving $(P_{f,X})$, an algorithm at every step can sequentially
call the oracle at a query point $x\in \bE$, obtaining the value $\cO(f,x)$.
In the sequel, we always assume the oracle to be {\sl local}, meaning
that for all $x\in \bE$ and $f,g\in\cF$  such that $f(\cdot)=g(\cdot)$
in a neighbourhood of $x$, we have $\cO(f,x)=\cO(g,x)$.
The most common example of oracle is the first-order
oracle, which returns the value and a subgradient of $f$ at $x$. However, observe that when the subdifferential
is not a singleton not every such oracle
satisfies the local property, and we need to further restrict it to satisfy
locality.

A $T$-step algorithm $\cM$, utilizing oracle $\cO$, for the family
$\cP(\cF,X)$ is a procedure as follows.  As applied to a problem
$(P_{f,X})$ with $f\in \cF$, $\cM$ generates a sequence
$x_t=x_t(\cM,f)$, $1\leq t\leq T$
of {\sl search points} according to the recurrence
$$
x_t=X_t(\{x_\tau,\cO(f,x_\tau)\}_{\tau=1}^{t-1}),
$$
where the {\sl search rules} $X_t(\cdot)$ are deterministic functions of
their arguments; we can identify $\cM$ with the collection of these rules.
Thus, $x_1$ is specified by $\cM$ and is independent of $f$, and all
subsequent search points are deterministic functions of the preceding
search points and the information on $f$ provided by $\cO$ when queried at these
points. We treat $x_T=x_T(\cM,f)$ as the approximate
solution generated by the $T$-step solution method $\cM$ applied to
$(P_{f,X})$, and define the {\sl minimax risk} associated with the family
$\cP(\cF,X)$ and oracle $\cO$ as the function of $T$ defined by
$$
\Risk_{\cF,X,\cO}(T)=\inf\limits_{\cM}\sup\limits_{f\in \cF}\left[f(x_T(\cM,f))-\Opt(f)\right],
$$
where the right hand side infinum is taken over all $T$-step solution
algorithms $\cM$ utilizing oracle $\cO$ and such that $x_T(\cM,f)\in X$ for all $f\in \cF$. The inverse to the risk function
$$
\cC_{\cF,X,\cO}(\varepsilon)=\min\left\{k:\Risk_{\cF,X,\cO}(k)\leq\varepsilon\right\}
$$
for $\varepsilon>0$ is called the {\sl information-based (or oracle) complexity of
the family $\cP(\cF,X)$} with respect to oracle $\cO$.\\

\noindent{\bf Geometry and Smoothness:} Let $\bE$  be an $n$-dimensional Euclidean space, and $\|\cdot\|$
be a norm on $\bE$ (not necessarily the Euclidean one). Let, further,
$X$ be a nonempty closed and bounded convex set in $\bE$.
Given
a positive real $L$ and $\kappa\in (1,2]$, consider the family
$\cF_{\|\cdot\|}(\kappa,L)$ of all continuously differentiable convex
functions $f:\bE\to\bR$ which are $(\kappa,L)$-smooth w.r.t. $\|\cdot\|$, i.e. satisfy the relation
\begin{equation}\label{eqLip}
 \|\nabla f(x)-\nabla f(y)\|_*\leq L\|x-y\|^{\kappa-1} \quad \forall x,y\in E,
\end{equation}
where $\|\cdot\|_*$ is the norm conjugate to $\|\cdot\|$. We  associate with
$\|\cdot\|,X,\kappa,L$ the family of convex optimization problems
$\cP=\cP(\cF_{\|\cdot\|}(\kappa,L),X)$.

We assume the family $\cF_{\|\cdot\|}(\kappa,L)$ is equipped with a local oracle $\cO$. To avoid extra words, we assume that this oracle is at least as powerful
as the First Order oracle, meaning that $f(x),\nabla f(x)$ is a component of $\cO(f,x)$.\par
Our goal is to establish lower bounds on the risk $\Risk(T)$, taken w.r.t.
the oracle $\cO$,  of the just defined family of problems $\cP$.
In the sequel, we focus solely on the {\sl`large-scale' case} $n\geq T$,
and the reason is as follows: it is known \cite{Nemirovski:1983} that
when $T\gg n$, $\Risk(T)$ ``basically forgets the details specifying
$\cP$'' and is upper-bounded by $O(\exp\{-CT/n\})$, where $C$
is an absolute constant, with the data $\|\cdot\|,X,L,\kappa$ 
of $\cP$ affecting only the hidden factor in the outer $O(\cdot)$
and thus
irrelevant when $T\gg n$. In contrast to this, in the large-scale regime $T\leq n$,
$\Risk(T)$ is (at least in the cases we are about to
consider) nearly independent of $n$ and goes to 0 {\sl sublinearly}
as $T$ grows, and its behavior in this range heavily depends on $\cP$.
In what follows, we focus solely on the large-scale regime.

\section{Local Smoothing} \label{sec:Smoothing}

In this section we introduce the main component of our technique, a Moreau-type  approximation of a
nonsmooth convex function $f$ by a smooth one. The main feature of this smoothing, instrumental for our ultimate goals, is
that it is local -- the local behaviour of the approximation at a point depends solely on the restriction of $f$ onto a neighbourhood of the point, the size of the neighbourhood being under our full control.

\subsection{Smoothing Kernel}\label{sect1}
Let $(\bE,\langle\cdot,\cdot\rangle)$ be a finite-dimensional Euclidean space,
$\|\cdot\|$ be a norm on $\bE$ (not necessarily induced by $\langle\cdot,\cdot\rangle$),
and $\cC_{\|\cdot\|}$ be the set of
all Lipschitz continuous, with constant 1 w.r.t. $\|\cdot\|$, convex functions on $\bE$. Let also
$\phi(\cdot)$ (``smoothing kernel'') be a twice continuously differentiable  convex function defined on an open convex
set $\Dom\phi \subset \bE$ with the following properties:
\begin{enumerate}
\item[A.]\label{item1} $0\in\Dom \phi$ and $\phi(0)=0$, $\phi'(0)=0$;
\item[B.]\label{item2} There exists a compact convex set $G\subseteq \Dom\phi$ such that $0\in\inter G$ and $\phi(x)>\|x\|$ for all $x\in\partial G$.
\item[C.] For some $M_\phi<\infty$ we have
\begin{equation}\label{Mphi}
\langle e,\nabla^2\phi(h) e \rangle \leq M_\phi \|e\|^2 \quad\forall (e\in \bE,h\in G).
\end{equation}
\end{enumerate}

Note that A and B imply that for all $f\in\cC_{\|\cdot\|}$, the function $f(x)+\phi(x)$
attains its minimum on the set $\inter G$. Indeed, for  every $x\in\partial G$ we
have $f(x)+\phi(x)\geq f(0)-\|x\|+\phi(x)>f(0)+\phi(0)$, so that the (clearly existing)
minimizer of $f+\phi$ on $G$ is a point from $\inter G$. As a result, for every
$f\in\cC_{\|\cdot\|}$ and $x\in \bE$ one has
\begin{equation}\label{eq256}
\min_{h\in\hbox{\scriptsize\rm Dom}\phi}[f(x+h)+\phi(h)]=\min_{h\in\hbox{\scriptsize\rm int} G} [f(x+h)+\phi(h)],
\end{equation}
and the right hand side minimum is achieved.

Given a function $f\in\cC$, we refer to the function
$$
\cS[f](x)=\min_{h\in\hbox{\scriptsize\rm Dom}\,\phi} [f(x+h)+\phi(h)]=\min_{h\in G}[f(x+h)+\phi(h)]
$$
as to the {\sl smoothing} of $f$. Observe that by our assumptions on $\phi$ we have
\begin{enumerate}
\item $\cS[f](x)=f(x+h(x))+\phi(h(x))$, where $h(x)\in \inter G$ is such that
\begin{equation}\label{eq1}
f'(x+h(x))+\phi'(h(x))=0
\end{equation}
for properly selected $f'(x+h(x))\in \partial f(x+h(x))$;
\item $f(x)\geq \cS[f](x)\geq f(x)-\rho_{\|\cdot\|}(G)$, where
$$
\rho_{\|\cdot\|}(G)=\max_{h\in G}\|h\|;
$$
indeed, by A we have $\phi(h)\geq \phi(0)=0$, so that $f(x)= f(x)+\phi(0)\geq \cS[f](x)=f(x+h(x))+\phi(h(x))\geq f(x+h(x)) \geq f(x)-\|h(x)\|$ (recall that $f\in\cC_{\|\cdot\|}$), while $h(x)\in G$.
\item We have that for all $f\in \cC$
\begin{equation}\label{then}
\|\nabla \cS[f](x)-\nabla\cS[f](y)\|_*\leq M_{\phi}\|x-y\|\quad\forall x,y\in E.
\end{equation}
\end{enumerate}
For a proof of (\ref{then}) see Section \ref{App1} in the Appendix.

\subsection{Approximating a function by smoothing}\label{sect2}
For $\chi>0$ and $f\in \cC_{\|\cdot\|}$, let
$$
\cS_\chi[f](x)=\min_{h\in\chi \hbox{\scriptsize\rm Dom}\,\phi} [f(x)+\chi\phi(h/\chi)].
$$
Observe that $\cS[f]_\chi(\cdot)$ can be obtained as follows:
\begin{itemize}
\item We associate with $f\in \cC_{\|\cdot\|}$ the function $f_{\chi}(x)=\chi^{-1}f(\chi x)$;
observe that this function belongs to $\cC_{\|\cdot\|}$ along with $f$;
\item We pass from $f_{\chi}$ to its smoothing
$$
\begin{array}{rcl}
\cS[f_{\chi}](x) &=& \min_{g\in\scriptsize{\Dom}\phi}\left[f_{\chi}(x+g)+\phi(g)\right] \\
  &=&\min_{g\in\scriptsize{\Dom}\phi}\left[\chi^{-1}f(\chi x+\chi g)+\phi(g)\right]\\
  &=&\chi^{-1}\min_{h\in\chi \scriptsize{\Dom}\phi}\left[f(\chi x+h)+\chi \phi(h/\chi)\right]\\
  &=&\chi^{-1}\cS_\chi[f](\chi x). \\
\end{array}
$$
\end{itemize}
It follows that
$$
\cS_\chi[f](x)=\chi\cS[f_{\chi}](\chi^{-1} x).
$$
The latter relation combines with (\ref{then}) to imply that
$$
\|\nabla\cS_\chi[f](x)-\nabla\cS_\chi[f](y)\|_*\leq \chi^{-1}M_{\phi}\|x-y\|\,\,\forall x,y.
$$

As bottom-line, if we can find a function $\phi$ as described above we have that for
any convex function $f:\bE\to\bR$ with Lipschitz constant 1 w.r.t. $\|\cdot\|$ and every $\chi>0$
there exists a smooth (i.e., with Lipschitz continuous gradient) approximation $\cS_\chi[f]$ that satisfies:

\begin{enumerate}
\item[S.1.] $\cS_\chi[f]$ is convex and Lipschitz continuous with constant 1 w.r.t. $\|\cdot\|$  and has a Lipschitz continuous gradient, with constant $M_{\phi}/\chi$, w.r.t. $\|\cdot\|$:
$$
\|\nabla\cS_\chi[f](x)-\nabla\cS[f](y)\|_*\leq \chi^{-1}M_{\phi}\|x-y\|\,\,\forall x,y;
$$
\item[S.2.] $\sup_{x\in E}|f(x)-\cS_\chi[f](x)|\leq \chi \rho_{\|\cdot\|}(G)$. Moreover,
$f(x)\geq \cS_\chi[f](x)\geq f(x)-\chi \rho_{\|\cdot\|}(G)$.
\item[S.3.] $\cS_\chi[f]$ depends on $f$ in a local fashion: the value and the derivative of
$\cS_\chi[f]$ at $x$ depends only on the restriction of $f$ onto the set $x+\chi G$.
\end{enumerate}

\subsection{Example: $p$-norm smoothing}\label{sectpnorm}
Let $n>1$ and $p\in[2,\infty]$, and consider the case of $\bE=\bR^n$,
endowed with the standard inner product, and  $\|\cdot\|=\|\cdot\|_p$.
Assume for a moment that $p>2$, and
let $r$ be a real such that $2< r\leq p$.  Let us select $\theta>1$ such that $2\theta/r<1$ and set
\begin{equation}\label{let11}
\begin{array}{rcl}
\phi(x)&=&\phi_{r}(x)=2\left(\sum_{j=1}^n|x_j|^r\right)^{2\theta/r},\\
G&=&\{x\in \bR^n:\|x\|_p\leq1\}.\\
\end{array}
\end{equation}

Observe that $\phi$ is twice continuously differentiable on $\Dom\phi=\bR^n$
function satisfying A. Besides this, $r\leq p$ ensures that $\sum_j|x_j|^r\geq1$
whenever $\|x\|_p=1$, so that $\phi(x)>\|x\|_p$ when $x\in\partial G$, which
implies B. Besides, by choosing $r=\min[p, 3\ln n]$ and $\theta>1$ close enough to 1,
C is satisfied for $M_{\phi} = O(1) \min[p,\ln n]$ (for a proof we refer to Section
\ref{App2} in the Appendix).

For the case of $p=2$, we can set $\phi(x)=2\|x\|_2^2$ and,
as above, $G=\{x:\|x\|_2\leq1\}$, clearly ensuring A, B, and the
validity of $C$ with $M_\phi=1$.

Applying the results of the previous section, we get

\begin{propos} \label{corinfty}
Let $p\in[2,\infty]$ and
$f:\bR^n\to\bR$ be a Lipschitz continuous, with constant
$1$ w.r.t. the norm $\|\cdot\|_p$, convex function. For
every $\chi>0$, there exists a convex continuously
differentiable function $\cS_\chi[f](x):\bR^n\to\bR$ with
the following properties:
\par
{\rm (i)} $f(x)\geq \cS_\chi[f](x)\geq f(x)-\chi$, for all $x$;
\par
{\rm (ii)} $\|\nabla \cS_\chi[f](x)-\nabla\cS_\chi[f](y)\|_{{p\over p-1}}\leq O(1)\min[p,\ln n]\chi^{-1}\|x-y\|_p$ for all $x,y$;
\par
{\rm (iii)} For every $x$, the restriction of $\cS_\chi[f](\cdot)$ on a small enough neighbourhood of $x$ depends solely on the restriction of $f$ on the set
$$
B^p_\chi(x)=\{y:\|y-x\|_p\leq \chi\}.
$$
\end{propos}


\section{Lower complexity Bounds for Smooth Convex Minimization} \label{sec:LCB}

In this section we utilize Proposition \ref{corinfty}  to prove our main result,
namely, a general lower bound on the oracle complexity of smooth
convex minimization, and then specify this result for the case of
minimization over $\|\cdot\|_p$ balls, where $2\leq p\leq\infty$.

\begin{propos} \label{main}
Let
\begin{enumerate}
\item[{\rm I.}] $\|\cdot\|$ be a norm on $\bR^n$ and $X$ be a nonempty
convex set containing the unit ball of $(\bR^n,\|\cdot\|)$;
\item[{\rm II.}] $T$ be a positive integer and $\Delta$ be a positive real
with the following property:
\\
One can point out $T$ linear forms $\langle\omega_i,\cdot\rangle$ on
$\bR^n$, $1\leq i\leq T$, such that
 \par (a) $\|\omega_i\|_*\leq 1$ for $i\leq T$, and
 \par
 (b)  for every collection $\xi^T=(\xi_1,...,\xi_T)$ with $\xi_i\in\{-1,1\}$, it holds
\begin{equation}\label{eqass1}
\min_{x\in X} \max_{1\leq i\leq k}\xi_i\langle\omega_i,x\rangle \leq -\Delta;
\end{equation}
\item[{\rm III.}] $M$  and $\rho$ be positive reals such that for properly selected
convex twice continuously differentiable on an open convex set
$\Dom \phi\subset\bR^n$ function $\phi$ and a convex compact subset
$G\subset\Dom\phi$ the triple $(\phi,G,M_\phi=M)$ satisfies  properties A, B, C
from Section \ref{sect1}  and $\rho_{\|\cdot\|}(G)\leq \rho$.
\end{enumerate}
Then for every $L>0$, $\kappa\in(1,2]$, every local oracle $\cO$ and every
$T$-step method $\cM$ associated with this oracle there exists a problem
$(P_{f,X})$ with $f\in \cF_{\|\cdot\|}(\kappa,L)$ such that
\begin{equation}\label{suchthat}
f(x_T(\cM,f)) -\thOpt(f)\geq {\Delta^\kappa\over 2^{\kappa+1}(\rho M)^{\kappa-1}}\cdot {L\over T^{\kappa-1}}.
\end{equation}
\end{propos}
{\bf Proof.} 1$^0$. Let us set
\begin{equation}\label{letusset}
\delta={\Delta\over 2T},\,\chi={\delta\over 2\rho}={\Delta\over 4T\rho},
\,\beta={L\chi^{\kappa-1}\over 2^{2-\kappa} M^{\kappa-1}}
={L\Delta^{\kappa-1}\over 2^\kappa(T\rho M)^{\kappa-1}}.
\end{equation}
\vspace{0.2cm}

2$^0$. Given a permutation $i\mapsto \sigma(i)$ of $\{1,...,T\}$ and
a collection $\xi^T\in\{-1,1\}^T$, we associate with these data the functions
$$
g^{\sigma(\cdot),\xi^T}(x)=\max_{1\leq i\leq T} \left[\xi_i\langle\omega_{\sigma(i)},x\rangle - (i-1)\delta\right].
$$
Observe that all these functions belong to $\cC_{\|\cdot\|}$ due to
$\|\omega_j\|_*\leq1$, for $j\leq T$, so that the smoothed functions
\begin{equation} \label{SmoothInstance}
f^{\sigma(\cdot),\xi^T}(x)=\beta \cS_\chi[g^{\sigma(\cdot),\xi^T}](x)
\end{equation}
(see Section \ref{sect2}) are well defined continuously differentiable  convex
functions on $\bR^n$ which, by item S.1 in Section \ref{sect2}, satisfy that
for all $x$, $y$ in $X$
$$
\|\nabla f^{\sigma(\cdot),\xi^T}(x)-\nabla f^{\sigma(\cdot),\xi^T}(y)\|_*\leq\beta\min[\chi^{-1}M\|x-y\|,2]
$$
whence
for all $x$, $y$ it holds
$$
\|\nabla f^{\sigma(\cdot),\xi^T}(x)-\nabla f^{\sigma(\cdot),\xi^T}(y)\|_{\ast}
\leq \beta 2^{2-\kappa}(\chi^{-1}M)^{\kappa-1}\|x-y\|^{\kappa-1}.
$$
Recalling the definition of $\beta$, we conclude that
$f^{\sigma(\cdot),\xi^T}(\cdot)\in\cF_{\|\cdot\|}(\kappa,L)$.

\vspace{0.2cm}

3$^0$. Given a local oracle ${\cal O}$ and an associated $T$-step
method $\cM$, let us define a sequence $x_1,\ldots,x_T$ of points
in $\bR^n$,
a permutation $\sigma(\cdot)$ of $\{1,...,T\}$ and a collection
$\xi^T\in\{-1,1\}^T$ by the following $T$-step recurrence:
\begin{itemize}
\item {\sl Step 1:}   $x_1$ is the first point of the trajectory of $\cM$
(this point depends solely on the method and is independent of the
problem the method is applied to). We define $\sigma(1)$ as the
index $i$, $1\leq i\leq T$, that maximizes $|\langle \omega_i,x_1\rangle|$,
and specify $\xi_1\in\{-1,1\}$ in such a way that
$\xi_1\langle \omega_{\sigma(1)},x_1\rangle = |\langle \omega_{\sigma(1)},x_1\rangle|$.
We set
$$
g^1(x)=\xi_1\langle \omega_{\sigma(1)},x\rangle,\, f^1(x)=\beta\cS_\chi[g^1](x).\\
$$
\item {\sl Step $t$, $2\leq t\leq T$:} At the beginning of this step,
we have at our disposal the already built points $x_\tau\in\bR^n$,
distinct from each other integers $\sigma(\tau)\in\{1,...,T\}$ and
quantities $\xi_\tau\in\{-1,1\}$, for $1\leq\tau <t$. At step $t$, we
build $x_t$, $\sigma(t)$, $\xi_t$, as follows.
We set
$$
g^{t-1}(x)=\max_{1\leq \tau<t} \left[\xi_\tau\langle\omega_{\sigma(\tau)},x\rangle - (\tau-1)\delta\right],
$$
thus getting a function from $\cC_{\|\cdot\|}$, and define its
smoothing $f^{t-1}(x)=\beta \cS_\chi[g^{t-1}](x)$ which,
same as above, belongs to $\cF_{\|\cdot\|}(\kappa,L)$. We further define
\begin{itemize}
\item $x_t$ as the $t$-th point of the trajectory of $\cM$ as applied to $f^{t-1}$,
\item $\sigma(t)$ as the index $i$ that maximizes $|\langle \omega_i,x_t\rangle|$,
over $i\leq T$ {\sl distinct from} $\sigma(1),...,\sigma(t-1)$,
\item $\xi_t\in\{-1,1\}$ such that $\xi_t\langle \omega_{\sigma(t)},x_t\rangle = |\langle \omega_{\sigma(t)},x_t\rangle|$
\end{itemize}
thus completing step $t$.
\end{itemize}
After $T$ steps of this recurrence, we get at our disposal a sequence
$x_1,\ldots,x_T$ of points from $\bR^n$, a permutation $\sigma(\cdot)$
of indexes $1,\ldots,T$ and a collection $\xi^T=(\xi_1,...,\xi_T)\in\{-1,1\}^T$;
these entities define the functions
$$
g^T=g^{\sigma(\cdot),\xi^T},\,f^T=\beta \cS_\chi[g^{\sigma(\cdot),\xi^T}].
$$


4$^0$. We claim that $x_1,\ldots,x_T$ is the trajectory of $\cM$ as
applied to $f^T$. By construction, $x_1$ indeed is the first point of the
trajectory of $\cM$ as applied to $f^T$. In view of this fact, taking into
account the definition of $x_t$ and the locality of the oracle $\cO$,
all we need  to support our claim is to verify that for every  $t$, $2\leq t\leq T$,
the functions $f^T$ and $f^{t-1}$ coincide in some neighbourhood of
$x_{t-1}$. By construction, we have that for $t\leq s \leq T$
\begin{equation}\label{maineq}
\xi_s\langle \omega_{\sigma(s)},x_{t-1}\rangle \leq |\langle\omega_{\sigma(t-1)},x_{t-1}\rangle|=\xi_{t-1}\langle \omega_{\sigma(t-1)},x_{t-1}\rangle,
\end{equation}
and thus
\begin{equation}\label{gk}
g^T(x)=\max\big[g^{t-1}(x),\underbrace{\max_{t\leq s\leq T}[\xi_s\langle \omega_\sigma(s),x\rangle -(s-1)\delta]}_{=:g_t(x)}\big]
\end{equation}
with
$$g^{t-1}(x_{t-1})\geq \xi_{t-1}\langle \omega_{\sigma(t-1)},x_{t-1}\rangle -(t-2)\delta.$$
Invoking (\ref{maineq}), we get
\begin{eqnarray*}
t\leq s\leq T &\Rightarrow& g^{t-1}(x_{t-1}) \geq [\xi_s\langle \omega_{\sigma(s)},x_{t-1}\rangle -(s-1)\delta]+\delta\\
			  &\Rightarrow&g^{t-1}(x_{t-1}) \geq g_t(x_{t-1})+\delta.
\end{eqnarray*}
Since both $g^{t-1}$ and $g_t$ belong to $\cC_{\|\cdot\|}$, it follows that
$g^{t-1}(x)\geq g_t(x)$ in the $\|\cdot\|$-ball $B$ of radius $\delta/2$
centered at $x_{t-1}$, whence, by (\ref{gk}),
$$
x\in B\quad\Rightarrow\quad g^T(x)=g^{t-1}(x).
$$
From $\chi \rho=\delta/2$ we have that $g^{t-1}\in\cC_{\|\cdot\|}$
and $g^T\in\cC_{\|\cdot\|}$ coincide on the set $x_{t-1}+\chi G$,
whence, as we know from item S.3 in Section \ref{sect2},
$f^{t-1}(\cdot)=\beta\cS_\chi[g^{t-1}](\cdot)$ and
$f^T(\cdot)=\beta\cS_\chi[g^T](\cdot)$ coincide in a neighbourhood of $x_{t-1}$, as claimed.

\vspace{0.2cm}

5$^0$. We have
\begin{eqnarray*}
g^T(x_T) &\geq& \xi_T\langle \omega_{\sigma(T)},x_T\rangle -(T-1)\delta \\
  &=& |\langle \omega_{\sigma(T)},x_T\rangle| -(T-1)\delta \\
  &\geq& -(T-1)\delta,
\end{eqnarray*}
whence, by item  S.2 in Section \ref{sect2}, $\cS_\chi[g^T](x_T)\geq -(T-1)\delta-\chi \rho \geq -T\delta=-\Delta/2$,
implying that
$$
f^T(x_T)\geq -\beta\Delta/2.
$$
On the other hand, by (\ref{eqass1}) there exists $x_{\ast}\in X$ such that
$g^T(x_{\ast})\leq \max_{1\leq i\leq T} \xi_i\langle\omega_{\sigma(i)},x_{\ast}\rangle\leq -\Delta$,
whence $\cS_\chi[g^T](x_{\ast})\leq g^T(x_{\ast})\leq-\Delta$ and thus
$\Opt(f^T)\leq f^T(x_{\ast})\leq-\beta\Delta$. Since, as we have seen,
$x_1,\ldots,x_T$ is the trajectory of $\cM$ as applied to $f^T$, $x_T$
is the approximate solution generated by $\cM$ as applied to
$f^T$, and we see that the inaccuracy of this solution, in terms of
the objective, is at least
${\beta\Delta\over 2}={\Delta^\kappa\over 2^{\kappa+1}(\rho M)^{\kappa-1}}\cdot {L\over T^{\kappa-1}},$
as required. Besides this, $f^T$ is of the form $f^{\sigma(\cdot),\xi^T}$,
and we have seen that all these functions belong to $\cF_{\|\cdot\|}(\kappa,L)$.  \qed

\begin{rem}
Note that the previous result immediately
implies the lower bound
$$ \thRisk_{\cF,X,\cO}(T) \geq
{\Delta^\kappa\over 2^{\kappa+1}(\rho M)^{\kappa-1}}\cdot {L\over T^{\kappa-1}},$$
on the complexity of the family $\cP(\cF_{\|\cdot\|}(\kappa,L),X)$, provided $X$ contains the unit $\|\cdot\|$-ball. Note that this bound is
independent of the local oracle $\cO$.
\par
The case where $X$ contains a $\|\cdot\|$-ball of radius $R>0$ instead of
the unit $\|\cdot\|$ ball can be reduced to the latter case by scaling instances
$f(\cdot)\mapsto f(\cdot /R)$, which corresponds to the transformation
$(\kappa,L)\mapsto(\kappa,\bar{L}:=LR^{\kappa})$ of the smoothness parameters.
 Thus, assuming that $X$ contains $\|\cdot\|$-ball of radius $R$, we have
\begin{equation}\label{RiskBound_Radius}
\thRisk_{\cF,R\cdot X,\cO}(T)\geq {\Delta^\kappa\over 2^{\kappa+1}(\rho M)^{\kappa-1}}\cdot {LR^{\kappa}\over T^{\kappa-1}}.
 \end{equation}
\end{rem}

\section{Case of $\|\cdot\|=\|\cdot\|_p$}
In this section we provide lower complexity bounds for
smooth convex optimization over $\|\cdot\|_p$-balls for the case when $\|\cdot\|=\|\cdot\|_p$.
In section \ref{FirstpBallcase} we show that Proposition \ref{main} implies
nearly tight optimal complexity bounds for the range
$2\leq p \leq \infty$; moreover, for fixed and finite $p$, the bound is tight within a factor depending solely on $p$. For the case
$p=\infty$, our lower bound matches the approximation
guarantees of the Conditional Gradient algorithm, up to a logarithmic
factor, proving near-optimality of the algorithm.

In section \ref{SecondpBallCase} we study the range $1\leq p <2$.
Here we prove nearly optimal complexity bounds by using
nearly-Euclidean sections of the $\|\cdot\|_p$-ball, together with
the $p=\infty$ lower bound.

\subsection{Smooth Convex Minimization over $\|\cdot\|_p$-balls, $2\leq p\leq \infty$}
\label{FirstpBallcase}
Consider the case when $\|\cdot\|$ is the norm $\|\cdot\|_p$
on $\bR^n$, $2\leq p\leq\infty$. Given positive integer $T\leq n$,
let us specify $\omega_i$, $1\leq i\leq T$, as the first $T$
standard basic orths, so that  for every collection
$\xi^T\in\{-1,1\}^T$ one clearly has
\begin{equation} \label{pNormEstimate}
\min_{\|x\|_p\leq 1}\max_{1\leq i\leq T} \xi_i\langle \omega_i,x\rangle \leq -T^{-1/p}.
\end{equation}
Invoking the results from Section \ref{sectpnorm} (cf. Proposition \ref{corinfty}),
we see that when $X\subset \bR^n$ is a convex set containing the unit
$\|\cdot\|_p$-ball, Assumptions II and III in Proposition \ref{main} are
satisfied with $M=O(1)\min[p,\ln n]$, $\rho=1$ and $\Delta=T^{-1/p}$.
Applying Proposition \ref{main}, we arrive at

\begin{cor}\label{cornormp}
Let $2\leq p \leq \infty$, $\kappa\in(1,2]$, $L>0$, and let $X\subset\bR^n$
be a convex set containing the unit ball w.r.t. $\|\cdot\|_p$.
Then, for every $T\leq n$ and every local oracle $\cO$,
the minimax risk of the
family of problems $\cP(\cF,X)$ with $\cF=\cF_{\|\cdot\|_p}(\kappa,L)$ admits
the lower bound
\begin{equation}\label{RiskBound}
\thRisk_{\cF,X,\cO}(T) = {\Omega(1)\over [\min[p,\ln n]]^{\kappa-1}} \dfrac{ L }{ T^{\kappa+{\kappa\over p}-1} },
\end{equation}
independent of the local oracle $\cO$ in use.
\end{cor}

Let us discuss some interesting consequences of the above result.

\par {\bf A. Complexity of smooth minimization over the box:} Corollary \ref{cornormp} implies that
when $X$ is the unit $\|\cdot\|_\infty$-ball in $\bR^n$, the $T$-step
minimax risk $\Risk_{\cF,X,\cO}(T)$ of minimizing over $X$ of objectives
from the family $\cF=\cF_{\|\cdot\|_\infty}(\kappa,L)$ in the range $T\leq n$
is lower-bounded by $\Omega(1/\ln n)L/T^{\kappa-1}$. On the other hand, from the standard efficiency estimate of Conditional Gradient
algorithm (see, e.g., \cite{Dem:Rub:1970,Pshe:1994,CJN:2013}) it follows that when applying the method to minimizing over $X$ a function $f\in \cF_{\|\cdot\|_\infty}(\kappa,L)$
over a convex compact domain $X$ of $\|\cdot\|_\infty$-diameter $2R$, the inaccuracy after $T=1,2,...$ steps does not exceed
$$
O(1) {LR^\kappa\over T^{\kappa-1}}
$$
We see that when $X$ is in-between two $\|\cdot\|_\infty$-balls with ratio of sizes $\theta$, the lower complexity bound coincides with the upper one within the
factor $O(1)\theta^\kappa\ln^{\kappa-1}(n)$. In particular, {\sl when minimizing functions $f\in\cF_{\|\cdot\|_\infty}(\kappa,L)$
over $n$-dimensional unit box $X$, the performance of the
Conditional Gradient algorithm, as expressed by its minimax risk,
cannot be improved by more than $O(\ln^{\kappa-1}(n))$ factor, for any
local oracle in use.} In fact, the same conclusion remains
true when $\|\cdot\|_\infty$ and the unit box $X$ are replaced with $\|\cdot\|_p$
and the unit $\|\cdot\|_p$-ball with ``large'' $p$, specifically, $p\geq \Omega(1)\ln n$.

\par {\bf B. Tightness:} In fact, in the case of $2\leq p<\infty$
the lower complexity bounds for smooth convex minimization over
$\|\cdot\|_p$-balls established in Corollary \ref{cornormp}, are tight:
it is shown in \cite{Nemirovski:1985}, see also \cite[Section 2.3]{Elster:1993}
that a properly modified  Nesterov's algorithm $\cN$ for smooth convex
optimization via the first-order oracle, as applied to problems of minimizing
functions $f$ from $\cF_{\|\cdot\|_p}(\kappa,L)$ over the $n$-dimensional unit $\|\cdot\|_p$-ball $X$, for any number $T\geq 1$ of steps ensures that
$$
f(x_T(\cN,f)) -\min_{x\in X} f(x)\leq C(p){L\over  T^{\kappa+{\kappa\over p}-1}},
$$
with $C(p)$ depending solely on $p$, which is in full accordance
with (\ref{RiskBound}).

\subsection{Smooth Convex Minimization over $\|\cdot\|_p$-balls, $1\leq p\leq 2$}
\label{SecondpBallCase}

We have obtained lower complexity bounds for smooth
convex minimization over $\|\cdot\|_p$-balls, where $2\leq p\leq \infty$.
Now we consider the case $1\leq p < 2$. We will build nearly
tight bounds by reducing to the case of $p= \infty$.

\begin{propos}\label{newcor}
Let $1\leq p \leq 2$, $\kappa\in(1,2]$, $L>0$, and let $X\subset \bR^n$
be a convex set containing the unit $\|\cdot\|_p$-ball.
For properly selected absolute constant  $\alpha\in(0,1)$ and for every $T\leq \alpha n$,
the minimax risk of
the family of problems $\cP(\cF,X)$ with $\cF=\cF_{\|\cdot\|_p}(\kappa,L)$
admits the lower bound
\begin{equation}\label{RiskBound_p_small}
\thRisk_{\cF,X,\cO}(T) =
\Omega\left( \dfrac{ L }{ \ln^{\kappa-1}(T+1) T^{\frac{3\kappa}{2}-1} } \right),
\end{equation}
independent of the local oracle $\cO$ in use.
\end{propos}

{\bf Proof.} {\bf 1$^0$.} By Dvoretzky's Theorem for the $\|\cdot\|_p$-ball \cite[Theorem 4.15]{Pisier:1989},
there exists an absolute constant $\alpha\in(0,1)$, such that for any positive integer $T\leq \alpha n$
there is a subspace $M\subseteq \bR^n$ of dimension $T$, and a centered
at the origin ellipsoid $E\subseteq M$, such that
\begin{equation} \label{incl}
{1\over 2}E \subseteq B_M:=\{x\in M:\|x\|_p\leq 1\} \subseteq E.
\end{equation}
Let $\{\gamma_i(\cdot):\, i=1,\ldots, T\}$ be linear forms on $M$ such that
$E=\{y\in M:\, \sum_{i=1}^{T} \gamma_i^2(y)\leq 1\}$. By the second inclusion
in (\ref{incl}), for every $i$, the maximum of the linear form $\gamma_i(\cdot)$
over $B_M$ does not exceed $1$, whence, by the Hahn-Banach Theorem,
the form $\gamma_i(\cdot)$ can be extended from $M$ to a linear form on the
entire $\bR^n$ to have the maximum over $B:=\{x:\|x\|_p\leq1\}$ not exceeding $1$. In other
words, we can point out vectors $g_i\in\bR^n$, $1\leq i\leq T$ such that
$\gamma_i(y)=\langle g_i, y\rangle$ for every $y\in M$ and
$\|g_i\|_{{p\over p-1}}\leq 1$, for all $1\leq i\leq T$.  Now consider the linear mapping
$$
x\mapsto Gx:=[\langle g_1,x\rangle;\ldots;\langle g_T,x\rangle]:\bR^n\to\bR^T.
$$
By the above, the operator norm of this mapping induced by the norms
$\|\cdot\|_p$ on the argument and $\|\cdot\|_\infty$ on the image spaces
does not exceed 1. As a result, when $f:\bR^T\to\bR$ belongs to
$\cF^T_{\|\cdot\|_\infty}(\kappa,L)$, the function $f_+:\bR^n\to\bR$
defined by $f_+(x)=f(Gx)$, for $x\in\bR^n$, belongs to the family
$\cF^n_{\|\cdot\|_p}(\kappa,L)$ \footnote{To avoid abuse of notation, we have
added to our usual notation $\cF_{\|\cdot\|}(\cdot,\cdot)$ for families of smooth
convex functions superscript indicating the argument dimension of the functions
in question.}. Setting $Y=GX$, we get a convex compact set in $\bR^T$.
\par
{\bf 2$^0$.}Observe that an optimization problem of the form
$$
\min_{y\in Y} f(y)\eqno{(P_{f,Y})}
$$
can be naturally reduced to the problem
$$
\min_{x\in X} f_+(x),\eqno{(P_{f_+,X})}
$$
and when the objective of the former problem belongs to
$\cF^T:=\cF^T_{\|\cdot\|_\infty}(\kappa,L)$, the objective of the latter problem
belongs to $\cF^n=\cF^n_{\|\cdot\|_p}(\kappa,L)$. It is intuitively clear that
the outlined reducibility implies that the complexity of solving problems from
the family $\Phi_n:=\{(P_{f,X}):f\in\cF^n\}$ cannot be smaller than the complexity
of solving problems from the family $\Phi_T:=\{(P_{f,Y}): f\in \cF^T\}$.
Taking this claim for granted (for a proof, see section \ref{App3}), let us derive
from it the desired result. To this end, observe that from the first inclusion in
(\ref{incl}) it follows that $Y$ contains the centered at the origin $\|\cdot\|_\infty$-ball
of radius $R={1\over2\sqrt{T}}$ (indeed, by construction this ball is already
contained in the image of ${1\over2}E\subset X$). By Corollary \ref{cornormp}
as applied to $p=\infty$ and to $Y$  in the role of $X$,
the worst-case, w.r.t. problems from the family $\Phi_T$, inaccuracy of
any $T$-step method based on a local oracle is at least \\
\begin{equation}
\label{Omega}{\Omega(1)\over \ln^{\kappa-1}(T+1)}{R^\kappa L\over  T^{\kappa-1}}
=\Omega\left({1\over \ln^{\kappa-1}(T+1)}{L\over T^{{3\kappa\over 2}-1}}\right),
\end{equation}
see (\ref{RiskBound_Radius}).
According to our claim, the latter quantity lower-bounds the worst-case,
w.r.t. problems from the family $\Phi_n$, inaccuracy of any $T$-step
method based on a local oracle, and (\ref{RiskBound_p_small}) follows. \qed
\par
Finally, we remark that the lower complexity bound stated in Proposition \ref{newcor}
in the smooth case $\kappa>1$ is, to the best of our knowledge, new
(the nonsmooth case $\kappa=1$ was considered already in \cite{Nemirovski:1983}).
This lower bound matches, up to logarithmic in $n$ factors, the upper complexity
bound for the family in question, see \cite{Elster:1993}.

\subsection{Matrix case}
We have proved lower bounds for smooth optimization over $\|\cdot\|_p$-balls
for all $1\leq p \leq \infty$.
Now we show how these bounds can be used for proving lower complexity
bounds on smooth convex minimization over Schatten norm balls in the spaces of matrices. Recall that the Shatten $p$-norm $\|x\|_{\Sch,p}$ of an $n\times n$
matrix $x$ is, by definition the $p$-norm of the vector of singular
values of $x$. The problems we are interested in now are of the form
$$ \min_{x\in \bR^{n\times n}} \{ f(x):\,\, \|x\|_{\Sch,p} \leq 1\}. $$
where $f\in\cF_{\|\cdot\|_{\Sch,p}}(\kappa,L)$.
\par
Observe that
Corollary \ref{cornormp} remains true when replacing
in it the embedding space $\bE=\bR^n$ of $X$ with the space $\bE=\bR^{n\times n}$
of $n\times n$ matrices, the norm $\|\cdot\|_p$ on $\bR^n$ with the Schatten
norm $\|\cdot\|_{\Sch,p}$, and the requirement ``$X\subset \bR^n$ is a convex set
containing the unit ball of $\|\cdot\|_p$'' with the requirement
``$X\subset\bR^{n\times n}$ is a convex set containing the unit ball of $\|\cdot\|_{\Sch,p}$.'' This claim is an immediate consequence of the fact that when
restricting an $n\times n$ matrix onto its diagonal, we get a linear mapping
of $\bR^{n\times n}$ onto $\bR^n$, and the factor norm on $\bR^n$ induced,
via this mapping, by $\|\cdot\|_{\Sch,p}$ is nothing but the usual $\|\cdot\|_p$-norm.
Consequently,  minimizing a function from $\cF_{\|\cdot\|_p}(\kappa,L)$ over
the unit $\|\cdot\|_p$ ball $X$ of $\bR^n$ reduces to minimizing a convex
function of exactly the same smoothness, as measured w.r.t. $\|\cdot\|_{\Sch,p}$,
over the unit Schatten $p$-norm ball $X^+$ of $\bR^{n\times n}$. As a result,
every universal (i.e., valid for every local oracle) lower bound on the
minimax risk for the problem class $\cP(\cF_{\|\cdot\|_p}(\kappa,L),X)$
automatically is a universal  lower bound on the minimax risk for the problem
class $\cP(\cF_{\|\cdot\|_{\Sch,p}}(\kappa,L),X^+)$.

Note, however, that the ``matrix extension'' of our lower complexity bounds is not ``completely costless'' --
the resulting bounds are applicable when $T\leq n$ ($p\geq2$) of $T\leq O(1)n$ ($1\leq p\leq 2$), and $n$ now is the square root of the actual dimension of $x$.
Thus, in the matrix case our lower complexity bounds are applicable in relatively  more narrow range of values of $T$.

\bibliographystyle{abbrv}
\bibliography{mybibfile}  

\appendix

\section{Appendix}
\subsection{Justification of (\ref{then})}\label{App1}

In order to prove (\ref{then}),
 by the standard approximation argument, it suffices to establish this relation
in the case when, in addition to the inclusion
$f\in \cC_{\|\cdot\|}$ and the assumptions A -- C on $\phi$,
 $f$ and $\phi$ are C$^\infty$ smooth and $\phi$ is strongly convex.
By (\ref{eq1}),
\begin{equation}\label{eq16}
\cS[f](x)=f(x+h(x))+\phi(h(x)),
\end{equation}
where $h:E\to G$ is well defined and solves the nonlinear system of equations
\begin{equation}\label{eq3}
F(x,h(x))=0,\,\,
F(x,h):= f'(x+h)+\phi'(h).
\end{equation}
We have ${\partial F(x,h)\over\partial h}=f''(x+h)+\phi''(h)\succ0$, implying
by the Implicit Function Theorem that $h(x)$ is smooth. Differentiating
the identity $F(x,h(x))\equiv 0$, we get
\begin{eqnarray*}
\underbrace{f''(x+h(x))}_{P}[I+h'(x)]+\underbrace{\phi''(h(x))}_{Q}h'(x)=0\\
\Leftrightarrow \quad P+(P+Q)h'(x)=0\\
\\
\Rightarrow h'(x)=-[P+Q]^{-1}P=[P+Q]^{-1}Q-I.\\
\end{eqnarray*}

On the other hand, differentiating (\ref{eq16}), we get
\begin{eqnarray*}
\langle \nabla \cS[f](x),e\rangle &=&\langle f'(x+h(x)),e+h'(x)e\rangle +\langle \phi'(h(x)),h'(x)e\rangle\\
 &=&\langle f'(x+h(x)),e\rangle +\langle \underbrace{f'(x+h(x))+\phi'(h(x))}_{=0},h'(x)e\rangle \\
 &=&-\langle \phi'(h(x)),e\rangle,
\end{eqnarray*}
that is,
$$
\nabla\cS[f](x) = -\phi'(h(x)).
$$
As a result, for all $e$, $x$, we have, taking into account that $P$, $Q$ are symmetric positive definite,
\begin{eqnarray*}
\langle e,\nabla^2\cS[f](x)e\rangle &=& -\langle h'(x)e,\phi''(h(x))e\rangle \\
 &=& -\langle [[P+Q]^{-1}Q-I]e, Qe \rangle\\
 &=& \langle e,Qe\rangle-\langle e,Q[P+Q]^{-1}Qe\rangle \\
 &\leq& \langle e,Qe \rangle \leq M_{\phi}\|e\|^2,
\end{eqnarray*}
and (\ref{then}) follows.

\subsection{Proof for section \ref{sectpnorm}}\label{App2}
When $x\neq0$, we have:
\begin{eqnarray}
\langle e, [\nabla^2\phi(x)]e \rangle &=&4r\theta (2\theta/r-1)(\sum_j|x_j|^r)^{2\theta/r-2}\left[\sum_j|x_j|^{r-1}\sign(x_j)e_j\right]^2 \nonumber\\
&& +4\theta(r-1)(\sum_j|x_j|^r)^{2\theta/r-1}\sum_j|x_j|^{r-2}e_j^2 \nonumber\\
\nonumber \\
&\leq&4\theta(r-1)(\sum_j|x_j|^r)^{2\theta/r-1}\sum_j|x_j|^{r-2}e_j^2 \label{thetaBound}\\
&\leq& 4\theta(r-1)\|x\|_r^{2\theta-r}\left[\sum_j|x_j|^{(r-2){r\over r-2}}\right]^{{r-2\over r}}\left[\sum_j|e_j|^{2{r\over2}}\right]^{{2\over r}} \label{nBound1}\\
&=&4\theta(r-1)\|x\|_r^{2\theta-2}\|e\|_r^2,\nonumber\\
\Rightarrow \langle e, [\nabla^2\phi(x)]e \rangle&\leq&4\theta(r-1)\|x\|_r^{2\theta-2}\|e\|_r^2\label{nBound2}\\
\nonumber
\end{eqnarray}
(we used that $2\theta/r<1$ in (\ref{thetaBound}) and the H\"older inequality in (\ref{nBound1})). By continuity, the resulting inequality holds true when $x=0$ as well.
\par
Now let us set $r=\min[p,3 \ln n]$. When $p\leq 3\ln n$, we have $r=p$, and (\ref{nBound2}) reads
$$
x\in G\quad\Rightarrow\quad \langle e,[\nabla^2\phi(x)]e \rangle \leq4\theta p\|e\|_p^2\,\,\,\forall e\in\bR^n,
$$
expressing the fact that $\phi$, $G$ satisfy assumption C with $M_\phi=5p$, provided that $1<\theta\leq 5/4$. When $p>3\ln n$, we have $r=3\ln n$, whence $\|e\|_r\leq n^{{1\over r}-{1\over p}}\|e\|_p \leq \exp\{1/3\}\|e\|_p$, so that for $\theta>1$ close enough to 1 (what is ``close enough'', depends solely on $n$) (\ref{nBound2}) reads
$$
x\in G\quad\Rightarrow\quad \langle e,[\nabla^2\phi(x)]e \rangle \leq 5r\|e\|_r^2\leq 15\exp\{2/3\}\ln(n)\|e\|_p^2\,\,\,\forall e\in\bR^n,
$$
expressing the fact that $\phi$, $G$ satisfy assumption C with $M_\phi=15\exp\{2/3\}\ln(n)$. Thus, $\phi$, $G$ satisfy assumption C with $M_\phi=O(1)\min[p,\ln(n)]$.  

\subsection{Item 2$^0$ of the proof of Proposition \ref{newcor}}\label{App3}
 Observe, first, that the claim we intend to justify indeed needs a justification: we cannot just argue that solving  ``lifted'' problems -- those from the family $\Phi_n^+=\{(P_{f^+,X}):f\in\cF^T_{\|\cdot\|_{\infty}}(\kappa,L)\}\subset \Phi_n$ -- cannot be simpler than solving problems from $\Phi_T$ due to the fact that the problems from the latter family can be reduced to those from the former one; we should specify the local oracles associated with the families in question, and to ensure that ``lifting'' does {\sl not} simplify problems just because the oracle for the ``lifted'' family is more informative than the oracle for the original family. The justification here is as follows: observe that among local oracles for families of real-valued functions on $\bR^m$ there is the ``most informative'' one, let us call it {\sl maximal}; when queried about a function $f$ at  a point $x$, the maximal oracle returns the class $\widetilde{f}$ of $f$ w.r.t. the equivalence relation ``$f$ is equivalent to $g$ if and only if $f$ and $g$ coincide with each other in some (perhaps depending on $f$ and $g$) neighbourhood of $x$.'' Clearly the maximal oracle allows to mimic any other local oracle, so that for every family of problems, the lower complexity bounds valid for the maximal oracle are valid for any other local oracle. Now, it is easily seen that the maximal oracle for the family of functions $\cF^T$ induces the maximal oracle for the lifted family $\{f^+:f\in\cF^T\}$; with this in mind, it is immediately seen that any maximal-oracle-based $T$-step method ${\cal M}^+$ for solving problems from the family $\Phi_n^+$ induces a maximal-oracle-based $T$-step method ${\cal M}$ for solving problems from the family $\Phi^T$ in such a way that the trajectory $x_1,x_2,\ldots$ of ${\cal M}$ on a problem $(P_{f,Y})$ is linked to the trajectory $x_1^+,x_2^+,\ldots$ of ${\cal M}^+$ on $(P_{f_+,X})$ by the relation $x_t=Gx_t^+$. Consequently, when the maximal oracles are used, a lower bound on the $T$-step minimax risk of $\Phi_T$ automatically is a lower bound on the same quantity for the family $\Phi_n^+=\{(P_{f_+,X}): f\in\cF^T\}$, and therefore for the larger family $\Phi_n$. In particular, the quantity (\ref{Omega}), which by Corollary \ref{cornormp} lower-bounds the maximal-oracle-based $T$-step minimax risk when solving problems from $\Phi_T$, lower-bounds the similar quantity for $\Phi_n$, and thus - the $T$-step minimax risk of $\Phi_n$ taken w.r.t. any local oracle, as claimed.

\end{document}